\documentclass[preprint]{imsart}

\pdfoutput=1
\RequirePackage[OT1]{fontenc}
\RequirePackage{amsthm,amsmath,amssymb}
\RequirePackage{natbib}
\usepackage{graphicx, subfig}
\usepackage{multirow}
\RequirePackage[colorlinks,citecolor=blue,urlcolor=blue]{hyperref}
\usepackage[margin=1in]{geometry}

\def\T{{ \mathrm{\scriptscriptstyle T} }}
\def\talpha{\tilde{\alpha}}
\def\R{\mathbb{R}}

\newtheorem{theorem}{Theorem}
\newtheorem{condition}[theorem]{Condition}

\begin{document}

\begin{frontmatter}
\title{Testing for Serial Dependence in Binomial Time Series \textrm{I}: Parameter Driven Models}
\runtitle{Testing for Serial Dependence in Binomial Mixed Models}
\author{W.T.M. Dunsmuir\\School of Mathematics and Statistics, University of New South Wales, Sydney, Australia\\ w.dunsmuir@unsw.edu.au\\
\and \\
J.Y. He\\School of Mathematics and Statistics, University of New South Wales, Sydney, Australia\\
jieyi.he@unsw.edu.au}

\maketitle

\begin{abstract}
Binomial time series in which the logit of the probability of success is modelled as a linear function of observed regressors and a stationary latent Gaussian process are considered. Score tests are developed to first test for the existence of a latent process and, subsequent to that, evidence of serial dependence in that latent process. The test for the existence of a latent process is important because, if one is present, standard logistic regression methods will produce inconsistent estimates of the regression parameters. However the score test is non-standard and any serial dependence in the latent process will require consideration of nuisance parameters which cannot be estimated under the null hypothesis of no latent process. The paper describes how a supremum-type test can be applied. If a latent process is detected, consistent estimation of its variance and the regression parameters can be done using marginal estimation which is easily implemented using generalised linear mixed model methods. The test for serial dependence in a latent process does not involve nuisance parameters and is based on the covariances between residuals centered at functions of the latent process conditional on the observations. This requires numerical integration in order to compute the test statistic. Relevant asymptotic results are derived and confirmed using simulation evidence. Application to binary and binomial time series is made. For binary series in particular, a complication is that the variance of the latent process, even if present, can be estimated to be zero with a high probability.
\end{abstract}


\end{frontmatter}

\section{Introduction} \label{Sec: Introduction}

This paper focusses on the development of methods for detecting serial dependence in time series of binomial counts in which the logit of the probability of success at each time point is a linear function of regression variables and a latent autocorrelated process. Previous work \cite{davis2000autocorrelation} and \cite{davis2009negative} focussed on analogous set-ups where the observed responses are conditionally Poisson or negative binomial respectively. In these situations, use of ordinary generalized linear model (GLM) fitting leads to consistent and asymptotically normal estimates of the regression parameters (other than the intercept) even if there is a latent autocorrelated process in the log-mean process. These results allow residuals based on the GLM estimates to be used to construct estimates and tests of serial dependence along the lines of standard time series practice. However, as recently discussed in \cite{wu2014parameter} and \cite{DunsHe2016}, for conditionally binomial responses, GLM estimates are inconsistent. \cite{wu2014parameter} proposed a semi-parametric method to obtain consistent estimates while \cite{DunsHe2016} proposed the use of marginal estimation implemented easily via generalized linear mixed model (GLMM) methods. Both approaches yield consistent and asymptotically normal estimates of the regression parameters. However, neither addresses the question of whether or not there is serial dependence in the latent random process. This paper focusses on developing and studying score-type tests for the detection of a latent process in binomial time series regression and, if a latent process is present, whether it is autocorrelated or not. If autocorrelation is detected then more advanced computationally intensive methods may be used to jointly model the regression effects and the variance and serial dependence parameters of the latent process -- see \cite{davisDunsmuir2016DVTS} for a review.

From a practical point of view it would be useful to have a test for a serially dependent latent process that was easy to implement without fitting a particular model for the serial dependence and for which large sample theory can be established. If the test suggested that a latent process was not present then the practitioner could use the results and conclusions based on a standard generalized linear model fit for inference about regression effects. If the test suggested a latent process which is independent is present then generalized linear mixed modelling fitting could be used to obtain estimates of the regression effects and the variability of the latent process (with caveats to be discussed below for the binary case). If serial dependence was detected in the latent process then the extra effect required to fit an appropriate autocorrelation model for the latent process would be justified.

Assume $x_{t}$ is an observed $r$-dimensional vector of regressors available at time $t$ and $\alpha_t$ is an unobserved stationary Gaussian process with zero mean, variance $\tau$ and autocovariances specified in terms of $\tau$ and serial dependence parameters $\psi$. Let
\begin{equation}\label{eq: LinPredWt}%
W_{t}=x_{t}^{\T}\beta+ \alpha_{t}
\end{equation}
be the state variable. Let $\theta=(\beta,\tau,\psi)$ be the list of all parameters. Given the $W_t$, $Y_{t}$ are assumed to be independent with density
\begin{equation} \label{eq: EFDensity}%
f(y_{t}|W_t)=\exp\left\{ y_{t}W_{t}- m_{t}b(W_{t}) + c(y_{t})\right\}, \quad c(y_{t}) = \log \binom{m_t}{y_t}.
\end{equation}
While the approach to testing for the existence of a latent process and any serial dependence in it that we present in this paper can be extended to the Poisson and negative binomial response cases of 
\eqref{eq: EFDensity}, we concentrate on the case where $Y_t$ represents the number of successes in $m_t$ binomial trials conducted at time $t$. Then $f(y_{t}|W_t)$ is the binomial distribution with $b(W_t)=\log(1+\exp(W_t))$,  $\mu_t=E(Y_{t}|W_{t})= m_{t}\dot{b}(W_{t})$ and 
$\sigma_t^2= \mathrm{Var}(Y_{t}|W_{t})= m_{t}\ddot{b}(W_{t})$.
\[
\pi_t =\dot{b}(W_t) = \frac{e^{W_t}}{1+e^{W_t}}, \quad \sigma_{t}^{2}= m_t \ddot{b}(W_t)=m_{t}\pi_{t}(1-\pi_{t})
\]
where $\dot b$ and $\ddot b$ denote first and second derivatives with respect to the argument of $b$.

The above model is often described as being ``parameter driven" using terminology of \cite{cox1981statistical}. In a companion paper \cite{DunsHeODtest2016} the GLARMA and BARMA type ``observation driven'' models are considered in which $\alpha_t$ is replaced by $Z_t$, a process defined conditionally on past values of itself and ``innovations" in terms of derivations of past observed responses $Y_t$. This leads to quite different test procedures and theoretical considerations than those considered in this paper.

The presence of $\alpha_t$ in the state variable $W_t$ distorts the variance-mean relationship in the exponential family density and that, at least in the Poisson and negative binomial case, induces overdispersion in the observed responses. There is a substantial literature on testing for overdispersion of a general type in exponential family models. However most of these tests are not based on time series models of the types considered here.

In this paper we develop specific methods for testing the null hypothesis of no latent process based on the alternative being the parameter driven model defined above. Under this null hypothesis the serial dependence parameters $\psi$ are nuisance parameters. To deal with the nuisance parameter issue we will apply the supremum test developed in \cite{davies1987hypothesis}. This method has been widely used and shown to be effective in accommodating nuisance parameters -- see \cite{andrews1994optimal}, \cite{andrews1996testing}, \cite{fokianos2012interventions} and \cite{calvori2014testing} for example. Further to the issue of nuisance parameters, in parameter driven models the asymptotic distribution of model based tests such as likelihood ratio and Wald tests, under the null hypothesis of no serial dependence, are non-standard with an approximate distribution more complicated than the chi-squared distribution. More detail about similar non-standard tests can be found in \cite{moran1971maximum} and \cite{self1987asymptotic}. This is further complicated by the current lack of any asymptotic theory for the maximum likelihood estimators as is noted in \cite{davisDunsmuir2016DVTS}. For these reasons this paper will develop implementable and theory based score type tests for serial dependence in binomial time series parameter driven models.

The remainder of the paper is organized as follows: Section \ref{Sec: Estimation Review} reviews some large sample properties of GLM and marginal estimation for the parameter driven model. Section \ref{Sec: Score Tests General} studies the score vectors of parameter driven models and proposed the two step score-type test for the existence of latent process and serial dependence; Section \ref{Sec: Simulations} assesses the accuracy of the asymptotic distributions of the two step score-type test; Section \ref{Sec: Binary/Binom TS Examples} illustrates the applications to a number of actual binomial time series, and investigates the large sample properties of the proposed statistics using simple simulation experiments.

\section{Estimation of Parameter driven models} \label{Sec: Estimation Review}

We assume there is sample of observations $y_1,\ldots,y_n$ and associated observations on the regressors $x_1,\ldots,x_n$ available for inference about the regression parameters and the parameters of the latent process. For the purposes of deriving the required tests it will be convenient to let $\talpha_t = \alpha_t/\tau ^{1/2}$. Note that whenever $\alpha_t$ is a Gaussian autoregressive moving average process so is $\talpha_t$ and, letting $\talpha=(\talpha_1,\ldots,\talpha_n)^{\T}$ for which
\begin{equation}
\textrm{Cov}\left( \talpha\right) = R(\psi) \label{eq: Cov tilde alpha}
\end{equation}
where $R$ is a correlation matrix.
We can then rewrite the state equation \eqref{eq: LinPredWt} as
\begin{equation}\label{eq: StateEqnTildeAlpha}%
W_{t}(\beta,\tau) = x_{t} ^{\T}\beta + \tau ^{1/2} \tilde \alpha_t.
\end{equation}
The log-likelihood is $l_n(\theta)=\log L_n(\theta)$ where
\begin{equation} \label{eq: FullLik tildealpha}%
L_{n}(\theta) = \int_{\mathbb{R}^{n}}f(y,\tilde\alpha; \theta) d\tilde \alpha,
\end{equation}
\begin{equation}\label{eq: fullLiklihood}
f(y,\tilde\alpha; \theta)=\exp\left\{ \sum_{t=1}^n y_t W_t(\beta, \tau) - m_t b(W_t(\beta,\tau))+c(y_t)\right\}g(\tilde \alpha, R(\psi))
\end{equation}
and $g(\tilde \alpha, R(\psi))$ is the multivariate normal density with zero mean vector and covariance matrix $R(\psi)$. To date, general methods for finding the maximum likelihood estimators, $\hat \theta$ for $\theta$ are not readily available nor are the usual consistency and asymptotic normality properties for $\hat \theta$ available even for simple models -- see the review \cite{davisDunsmuir2016DVTS}.

For Poisson and negative binomial response distributions the regression parameters $\beta$ can be estimated consistently using generalized linear modelling -- see \cite{davis2000autocorrelation}, \cite{davis2009negative}. This is equivalent to ignoring the latent process and estimating $\beta$ by maximizing the GLM log-likelihood
\begin{equation}\label{eq: loglikglm}%
l_0(\beta) =\log L_0(\beta) = \sum_{t=1} ^{n}\left[y_{t}(x_{t} ^{\T}\beta) - m_{t}b(x_{t}^{\T}\beta) + c(y_{t})\right].
\end{equation}
We let $\hat \beta^{(0)}$ denote the value of $\beta$ which maximises \eqref{eq: loglikglm}. \cite{DunsHe2016} show that, for the binomial distribution, $\hat\beta^{(0)}$ converges to $\beta'$,
the unique vector that solves
\begin{equation}\label{eq: betaprime equation}%
\underset{n\to\infty}\lim n^{-1}\sum_{t=1}^n m_{t}\left[ \int_{\mathbb{R}} \dot{b}(x_{t}^{\T}\beta_0 + \alpha)g(\alpha;\tau_0)d\alpha - \dot{b}(x_{t}^{\T}\beta^\prime)\right] x_{t}=0
\end{equation}
and $\beta'\ne \beta_{0}$.

To overcome the inconsistency observed in GLM estimation, \cite{DunsHe2016} propose use of marginal likelihood estimation, which maximises the likelihood constructed under the assumptions that the process $\alpha_t$ consists of independent identically distributed random variables. Under this assumption the full likelihood \eqref{eq: fullLiklihood} is replaced by the ``marginal" likelihood, and
\begin{equation}\label{eq: marginal likelihood}%
l_1(\delta) = \sum_{t=1}^{n} \log \int_{\R} \exp\left( y_{t}W_t(\beta,\tau)-m_t b(W_t(\beta,\tau)) + c(y_{t})\right) g(\talpha_t) d\talpha_t
\end{equation}
where $\delta = (\beta,\tau)$ and $g(\cdot)$ is the standard normal distribution. Let $\hat\delta^{(1)}$ be the estimates of the true parameters $\delta_0 = (\beta_0, \tau_0)$ obtained by maximising \eqref{eq: marginal likelihood} over the compact parameter space $\Theta=\{\beta \in \R^r: \|\beta - \beta_0\|
\le d_1\}\bigcap \{\tau\ge 0: |\tau - \tau_0| \le d_2\}$, where $d_1<\infty$, $d_2<\infty$. Marginal likelihood estimators of $\hat\delta^{(1)}$ can be easily obtained with standard software packages for fitting generalized linear mixed models. However, it is also demonstrated in \cite{DunsHe2016} that marginal estimation can result in $\hat\tau^{(1)}=0$ with a substantial probability in finite samples. In particular, for binary response data, $P(\hat\tau^{(1)} =0)$ can be close to $50\%$. Explanations for why this can be expected to occur more frequently for binary responses are provided in \cite{DunsHe2016}. We return to this point when considering the simulations and analysis of real data in later sections.

\section{Score Tests} \label{Sec: Score Tests General}

\subsection{Score Vector of Parameter Driven Models}

In this section we will derive the score vector associated with the likelihood in \eqref{eq: FullLik tildealpha} with respect to the various components of $\theta$. In section \ref{Sec: ST LatPro} we derive a score type test for testing $H_0: \tau\ne 0$ (a latent process exists). This test involves nuisance parameters ($\psi$ describing serial dependence) which cannot be estimated under this null hypothesis. In section \ref{Sec: ST SerDep} we propose a score test for detecting serial dependence in the latent process given $\tau \ne 0$. The null hypothesis of no serial dependence is specified by $H_{0}:\psi=0$.

First consider the component of the score vector associated with $\beta$:
\begin{equation*}
S_{\beta}(\theta) =\frac{\partial l_n(\theta)}{\partial \beta} =\frac{  \int_{\mathbb{R}^{n}}  \sum_{t=1} ^{n}[y_{t}- m_{t} \dot{b}(W_{t}(\beta,\tau))] x_{t} f(y,\tilde \alpha) d \tilde \alpha} {L_n(\theta)}
\end{equation*}
Under $H_0: \tau=0$,
\begin{equation*}
S_{\beta}(\beta,0,\psi) =  \sum_{t=1} ^{n}e_t(\beta,0) x_{t}
\end{equation*}
and $e_t(\beta,\tau)=y_{t}- m_{t}\dot{b}(W_t(\beta,\tau))$ are unscaled residuals. Note that setting $S_{\beta}(\beta,0,\tau)$ to zero and solving for $\beta$ gives the GLM\ estimate $\hat \beta^{(0)}$, hence $S_\beta(\hat \beta^{(0)}, 0, \psi)=0$ regardless of the serial correlation parameter $\psi$.

Under $H_0: \tau=0$ both numerator and denominator of the derivative of the log-likelihood $l_n(\theta)$ with respect to $\tau$ are expressions which tend to $0$ as $\tau \rightarrow 0$ so application of L'H\^{o}pital's rule gives, after some algebra,
\begin{align}\label{eq: ScoreVec Ln}%
S_{\tau}(\beta,0,\psi) & =\frac{1}{2} \sum_{t=1}^n [e_t(\beta,0)^2 - m_t \ddot b(W_t(\beta,0))] + \frac{1}{2}\sum_{t=1}^n \sum_{s \ne t}^n e_t(\beta,0)e_s(\beta,0)R(s,t;\psi)\nonumber\\
&= S_{\tau,1}(\beta)+S_{\tau,2}(\beta,\psi)
\end{align}
The first component $S_{\tau,1}(\beta)$ depends only on the regression parameters and under the null hypothesis these can be estimated consistently using GLM. The second expression,
$S_{\tau,2}(\beta,\psi)$, depends on the correlation parameters $\psi$ which are not estimable under the null hypothesis, we propose the use of the supremum type tests discussed in \cite{davies1987hypothesis} to deal with these nuisance parameters. Note that $S_{\tau,1}(\beta)$ and $S_{\tau,2}(\beta)$ are uncorrelated.

To set up the second part of the test connected with a latent process, namely the test for serial dependence, note that if $\tau >0 $ but $H_0: \psi=0$ is true (that is, there is a latent process but it is serially uncorrelated) then $S_{\beta}(\beta, \tau, 0)$ and $S_{\tau}(\beta,\tau, 0)$ are the components of the score vector corresponding to the marginal likelihood in \cite{DunsHe2016}. Also, the component of the score vector associated with the typical element $\psi_a$ is
\begin{equation*}
S_{\psi_a}(\theta)=\frac{1}{L_n(\theta)} \int_{\mathbb{R}^{n}} \frac{1}{2}\textrm{tr}\left(\frac{\partial R(\psi)}{\partial \psi_a} [R(\psi)^{-1}\tilde \alpha \tilde \alpha^\top R(\psi)^{-1}-R(\psi)^{-1}]\right)f(y,\tilde\alpha; \theta) d \tilde \alpha
\end{equation*}
where $f(y,\tilde\alpha; \theta)$ is defined in \eqref{eq: fullLiklihood}.
Note that $S_{\beta}(\hat \beta^{(1)}, \hat \tau^{(1)},0) = 0$ and $S_{\tau}(\hat \beta^{(1)}, \hat \tau^{(1)},0) = 0$, also, because the diagonal elements of $R$ are all unity, $\left[\partial R(\psi)/\partial \psi_a \right]_{t,t}=0$ and $R(0)=I_n$. We have, hence, under the null hypothesis of independence of the $\tilde \alpha_t$, $\tilde\alpha_{t-a}|y_{t-a}$ and $\tilde\alpha_{t}|y_{t}$
are also independent so that
\begin{equation}\label{eq: ScoreVec Ln psi}%
S_{\psi_a}(\beta,\tau,0) =\frac{1}{2}\sum_{s=1}^n \sum_{t \ne s}^{n} E(\tilde \alpha_s|y_s) E(\tilde \alpha_t|y_t)\left[\frac{\partial R(\psi)}{\partial \psi_a} \right]_{s,t}.
\end{equation}

Further simplification can be given when $\tilde \alpha_t$ is an autoregressive moving average process -- see section \ref{Sec: ST SerDep}.

\subsection{Score Test for existence of latent process}\label{Sec: ST LatPro}%

Based on the score with respect to $\tau$ given by \eqref{eq: ScoreVec Ln} the statistic
\begin{equation}\label{eq: Q tau theoretical}
\hat Q_{\tau}(\psi) = \frac{n^{-1}S_{\tau}^2(\hat \beta^{(0)},0,\psi)}{n^{-1}\mathrm{Var}(S_\tau(\hat \beta^{(0)},0,\psi))}
\end{equation}
can be used to test the null hypothesis of no latent process. For this we need the distribution of $n^{-1/2}S_\tau(\hat\beta^{(0)},0,\psi)$. Some regularity conditions are required:

\begin{condition}\label{Cond: mt} The sequence of binomial trials $\{m_{t}:1\le m_{t}\le M\}$ is:
\begin{description}
 \item (a) stationary for which $\kappa_j= P(m_t=j)$, $\kappa_M>0$, $\sum_{j=1}^{M}\kappa_j = 1$ and is a strongly mixing process (as defined in \cite{gallant1988unified}) independent of $\{X_t\}$ with mixing coefficients satisfying $\sum_{h=0}^\infty \nu(h)<\infty$; or

 \item (b) asymptotically stationary and for which the $\kappa_j$ are the limits of finite sample frequencies of occurrences $m_t = j$.
\end{description}
\end{condition}
Both specifications cover the simple case of a constant number of trials at each time $t$, $m_t=M$.

\begin{condition}\label{Cond: Xt}
The regression sequence is specified in one of two ways
\begin{description}
		\item (a) \textit{Deterministic covariates} defined with functions:  $x_{nt}=h(t/n)$ for some specified piecewise continuous vector function $h: [0,1]\to \mathbb{R}^r$, or,
		
		\item (b) \textit{Stochastic covariates} which are a stationary vector process: $x_{nt}=x_t$ for all $n$ where $\{x_t\}$ is an observed trajectory of a stationary process for which $E(e^{s^T X_t}) <\infty$ for all $s\in \mathbb{R}^r$.
	\end{description}
\end{condition}

\begin{condition}\label{Cond: Reg Full Rank}
Let $r=\dim(\beta)$. The regressor space $\mathbb{X}=\{x_{nt}: 1\le t\le n\}$ has $\texttt{rank}(\texttt{span}(\mathbb{X}))=r$ for sufficiently large $n$.
\end{condition}

\begin{theorem} \label{Thm: Nulldist existence latent}%
Given Conditions \ref{Cond: mt} to \ref{Cond: Reg Full Rank}, under the null hypothesis $H_{0}:\tau=0$, for any fixed $\psi$ such that $\sum_{h=1}^\infty \vert R(h;\psi)\vert < \infty$, as $n\to \infty$,
$\hat Q_{\tau}(\psi)\rightarrow \chi^2(1)$ in distribution.
\end{theorem}

For any fixed value $\psi$, the asymptotic normality of $n^{-1/2}S_\tau(\hat\beta^{(0)},0,\psi)$ can be established based on the asymptotic normality of GLM\ estimates, $\hat\beta^{(0)}$. Proof that the GLM estimates are consistent and asymptotically normal under $H_0: \tau=0$ follows the functional limit theorem approach used in \cite{davis2000autocorrelation} and \cite{davis2009negative} for example. Details are available in \cite{DunsHe2016}.

Under the null $H_{0}: \tau=0$, let $\theta_0=(\beta_0,0,\psi)$. Write
\[
n^{-1/2}S_{\tau}^{\dag}(\theta_{0})= n^{-1/2}S_{\tau,1}(\beta_0) - J_n^T(\beta_0) I_n^{-1}(\beta_0) U_{n}(\beta_0)+ n^{-1/2}S_{\tau,2}(\beta_0,\psi)
\]
where $U_{n}(\beta_{0})=n^{-1/2}\sum_{t=1}^n e_{t}(\beta_{0},0)x_{t}$, and
\begin{equation*}
J_n(\beta_0)= -\frac{1}{2n}\sum_{t=1}^n m_tb^{(3)}(x_t^{\T}\beta_0)x_t, \quad
I_n(\beta_0) = \frac{1}{n}\sum_{t=1}^n m_t \ddot{b}(x_t^{\T}\beta_0)x_tx_t^{\T}.
\end{equation*}
Then the essential steps of the proof of Theorem \ref{Thm: Nulldist existence latent} are to show that, for any fixed $\psi$, \newline $n^{-1/2}\left[S_{\tau}(\hat\beta^{(0)},0,\psi) - S_{\tau}^{\dag}(\theta_0)\right] =o_{p}(1)$ and $n^{-1/2}S_{\tau}^{\dag}(\theta_0)\to N(0, \underset{n\to\infty}\lim \textrm{Var}(n^{-1/2}S_{\tau}^{\dag}(\theta_{0}))$ in distribution. An outline proof is provided in the appendix.

Let $V_n(\beta_0,\psi)=\textrm{Var}(n^{-1/2}S_{\tau}^{\dagger}(\theta_0))= V_{n,1}(\beta_0)+V_{n,2}(\beta_0,\psi)$, where $V_{n,1}(\beta_0)= K_n(\beta_0)- J_n^{\T}(\beta_0)I_n^{-1}(\beta_0)J_n(\beta_0)$,
\begin{equation*}
K_n(\beta_0)=\frac{1}{4n}\sum_{t=1}^n \sigma_t^2(\beta_0)(1+(2-6/m_t)\sigma^2_t(\beta_0)),
\end{equation*}
and $V_{n,2}(\beta_0,\psi)=n^{-1}\sum_{t=2}^n \sigma_t^2(\beta_0)\sum_{h=1}^{t-1} R^2(h;\psi)\sigma_{t-h}^2(\beta_0)$.

In practice we replace $n^{-1}\textrm{Var}(S_{\tau}(\hat\beta^{(0)},0,\psi))$ in
\eqref{eq: Q tau theoretical} by $V_n(\hat\beta_{0}, \psi)$ (evaluated using GLM estimates)
 to give
\begin{equation}\label{eq: Q tau empirical}
\hat Q_{\tau}(\psi) = n^{-1}S_{\tau}^2(\hat \beta^{(0)},0,\psi)/V_{n}(\hat \beta^{(0)},\psi).
\end{equation}
Since $\hat\beta^{(0)}-\beta_{0}=o_{p}(1)$, under Condition \ref{Cond: Xt}, $V_{n}(\hat\beta^{(0)},\psi)-V_{n}(\beta_{0},\psi)=o_{p}(1)$ can be established. Then,
as a corollary to Theorem \ref{Thm: Nulldist existence latent}, $\hat Q_{\tau}(\psi)$ in
\eqref{eq: Q tau empirical} also converges to $\chi^2(1)$ in distribution.

We note in passing that $K(\beta_0)$, $J(\beta_0)$ and $I(\beta_0)$, which are limits of $K_n(\beta_0)$, $J_n(\beta_0)$, $I_n(\beta_0)$ as $n\to\infty$ respectively, can be calculated as closed form integrals for both the deterministic and stochastic specification of regressors. In practice, $\beta_0$ is not known and so we substitute $\hat \beta^{(0)}$ for $\beta_0$ in the above definitions. It is straightforward to show that the resulting quantities also converge to the relevant limits $K(\beta_0)$, $J(\beta_0)$ and $I(\beta_0)$. Again, if asymptotic formula are used to estimate the variance $V_{n}(\beta_{0},\psi)$ the large sample chi-squared distribution will also hold. We do not use these asymptotic formulae in this paper and all results are based on using the finite sample formulae evaluated at $\hat \beta^{(0)}$.

There are various ways in which the score test of no latent process based on $\hat Q_{\tau}(\psi)$ can be implemented. The first is to simply assume that, under $H_0: \tau=0$, $\psi=0$ also and the test is against the alternative that the latent process is an independent and identically distributed sequence. Then
\begin{equation}\label{eq: std ST}%
\hat Q_{\tau}(0) = n^{-1}S_{\tau,1}(\hat\beta^{(0)})/V_{n,1}(\hat\beta^{(0)})
\end{equation}
since when $\psi=0$, $R(s,t; 0) =0$ if $s\ne t$ and $S_{\tau,2}(\hat\beta^{(0)},0) = 0$.
We call this the ``standard" score test of no latent process and note that this is the same as the test for homogeneity often proposed in the literature for overdispersion. However, as simulations reported below demonstrate, $S_{\tau,1}(\beta)$ may not have power increasing to 1 as $\tau$ grows, particularly for binary responses.

The second option is relevant when a particular direction and strength of serial dependence in the latent process is being tested against. For example, $R(s,t;\psi) = \psi^{\vert s-t\vert}$ corresponding to an \textsc{ar}$(1)$ process. In this type of situation, $R(s,t;\psi)$ is capable of reflecting the type and strength of serial dependence of a latent process if $H_0:\tau=0$ is rejected and $\psi$ would be fixed at a particular value. Then $\hat Q_{\tau}(\psi)$ is used in Theorem \ref{Thm: Nulldist existence latent}. We refer to this as the score test against a particular ``choice" of serial dependence. This type of test will have limited application since a test based on a fixed value of $\psi$ will not likely be powerful against others.

The third version is based on some particular functional with respect to $\psi$ of $Q_\tau(\psi)$. We illustrate the method of \cite{davies1987hypothesis} for the supremum score test
\begin{equation}\label{eq: sup ST}%
\hat Q_\tau(\Psi) = \sup_{\psi \in \Psi} \hat Q_{\tau}(\psi)
\end{equation}
where $\Psi$ is a compact set for $\psi$ chosen to give a reasonable range of alternatives.

\cite{davies1987hypothesis} derived an upper bound for the upper tail probability of the supremum score test with a single nuisance parameter $\psi$ over $\Psi=[\psi_{\mathcal{L}}, \psi_{\mathcal{U}}]$ as
\begin{equation} \label{eq: Sup chi2 AsyDist}%
P\left\{\underset{\psi_{\mathcal{L}} \le \psi \le \psi_{\mathcal{U}}}\sup S(\psi) > u \right\}
\le P(\chi^2_{s} > u) + \int_{\psi_{\mathcal{L}}}^{\psi_{\mathcal{U}}}\varphi(\psi)d\psi
\end{equation}
in which $S(\psi)= Z_{1}^2(\psi) + \ldots + Z_{s}^2(\psi)$, where $Z_{i}(\psi)\sim N(0,1)$ for all $i=1,\ldots,s$, and $\varphi(\psi)$ can be calculated by
\begin{equation*}
\varphi(\psi) = \frac{1}{\sqrt{2\pi}}\int_{0}^{\infty} \left\{ 1- \prod_{j=1}^{s} (1+\lambda_{j}(\psi)t)^{-1/2}\right\}t^{-3/2}dt u^{\frac{s-1}{2}}e^{-\frac{u}{2}} \pi^{-\frac{1}{2}}2^{-\frac{s}{2}}/\Gamma(\frac{s}{2}+\frac{1}{2})
\end{equation*}
where $\lambda_{j}(\psi)$, $j=1,\ldots,s$ are the eigenvalues of the matrix $B(\psi)- A^{T}(\psi)A(\psi)$. Here $Y(\psi) = \partial Z(\psi)/\partial \psi$, and
\begin{equation*}
\mathrm{Var}\binom{Z(\psi)}{Y(\psi)} =
\begin{bmatrix}
I & A(\psi)\\
A^{T}(\psi) & B(\psi)
\end{bmatrix}.
\end{equation*}
For the simple case, $s=1$,
\begin{equation} \label{eq: Prob upcrossing dim one}%
\int_{\psi_{\mathcal{L}}}^{\psi_{\mathcal{U}}}\varphi(\psi)d\psi = \pi^{-1} e^{-\frac{u}{2}} \int_{\psi_{\mathcal{L}}}^{\psi_{\mathcal{U}}}\lambda^{1/2}(\psi)d\psi.
\end{equation}

Next we focus on the specific example when $\alpha_t$ is an \textsc{ar}$(1)$ process. For any fixed $\psi$, based on Theorem \ref{Thm: Nulldist existence latent}, $Z(\psi)= n^{-1/2} S^{\dag}(\beta_0,0,\psi)V_{n}(\beta_0, \psi)^{-1/2}$ is asymptotically normal. Note also $\hat Q_{\tau}(\psi)$ converges to $Z^2(\psi)$ in distribution, and $\textrm{Cov}(Z(\psi), Y(\psi))=0$, thus the distribution of \eqref{eq: Sup chi2 AsyDist} can be rewritten as $P\left\{\sup_{\psi_{\mathcal{L}} \le \psi \le \psi_{\mathcal{U}}} \hat Q_{\tau}(\psi) > u \right\} = \mathcal{F}_{\Psi}(u)$ where
\begin{equation}\label{eq: Sup S1 dist}%
\mathcal{F}_{\Psi}(u)= P(\chi^2_{1} > u) + \pi^{-1}e^{-\frac{u}{2}} \int_{\psi_{\mathcal{L}}}^{\psi_{\mathcal{U}}} \lambda^{1/2}(\psi)d\psi%
\end{equation}
in which
\begin{equation*}
\lambda(\psi)= \frac{1}{n}\sum_{h=1}^{n-1}\sum_{t=1}^{n-h}\sigma_{t}^2(\beta_0) \sigma_{t+h}^2(\beta_0) h^2\psi^{2(h-1)}V_{n}(\beta_0, \psi)^{-1} - \left[\frac{1}{n}\sum_{h=1}^{n-1}\sum_{t=1}^{n-h} \sigma_{t}^2(\beta_0)\sigma_{t+h}^2(\beta_0)h\psi^{2h-1}V_{n}(\beta_0, \psi)^{-1}\right]^2
\end{equation*}
since the integral in \eqref{eq: Sup S1 dist} does not have a closed form it is evaluated using  numerical approximation.

\subsection{Score Test for Serial Dependence}\label{Sec: ST SerDep}%

When $\tau > 0$, the null hypothesis of no serial dependence is specified by $H_{0}^{(\psi)}: \psi=0$. To implement a score test based on \eqref{eq: ScoreVec Ln psi} we specialise to the case where $R_n(\psi)$ is the $n\times n$ dimensional auto-covariance matrix for an autoregressive-moving-average (ARMA) process $\tilde\alpha_t$, which has unit variance for the innovations and parameters $\psi=(\phi^{T},\theta^{T})^{T}$ which represent the autoregressive and moving average coefficients. Defining the autocovariances in terms of the spectral density for an ARMA process, differentiating with respect to $\theta_a$ or $\phi_a$ and evaluating at $\psi=0$ results in

\begin{equation}
\left[\frac{\partial R(\psi)}{\partial \psi_a}|_{\psi=0} \right]_{s,t}=
\begin{cases} 1, &  |s-t|=a \\
0, &  |s-t|\ne a.
\end{cases}
\end{equation}
so that \eqref{eq: ScoreVec Ln psi} simplifies to
\begin{align*}
S_{\psi_a}(\beta,\tau,0)
& = \sum_{t=a+1}^n E(\tilde\alpha_{t}|y_{t}; \beta,\tau) E(\tilde\alpha_{t-a}|y_{t-a};\beta,\tau)
\end{align*}
Using integration by parts it is straightforward to show that for the
exponential family with canonical link,
\[
E(\tilde\alpha_{t}|y_{t};\beta,\tau)= \sqrt{\tau} \left[ y_{t}- m_{t}E(\dot{b}(x_{t}^{\T}\beta
+ \sqrt{\tau} \tilde\alpha_{t})|y_{t};\beta, \tau)\right].
\]
Now let
\begin{equation}\label{eq: Ut}%
U_{t}=y_{t}- m_{t}E(\dot{b}(x_{t}^{\T}\beta + \sqrt{\tau} \tilde\alpha_{t})|y_{t};\beta, \tau)
\end{equation}
so that the score function is%
\[
S_{\psi_a}(\beta,\tau,0)=\tau \sum_{t=1}^{n}U_{t}U_{t-a}%
\]
Note that this is exactly the same for the $a$th lag moving average parameters and the $a$th lag autoregressive parameter. That is, moving average terms and autoregressive terms have the same contribution to the score and hence to a score statistic constructed from it. Hence a score statistic for detecting serial dependence from an \textsc{arma}$(p,q)$ model which has no common roots in the autoregressive and moving average polynomials can be constructed using $S_{\psi_a}(\beta,\tau, 0)$, $a=1,\ldots, L$ for some integer $L = \max(p,q)$. This is a pure significance test as described in \cite{poskitt1980testing} for the analogous score test of \textsc{arma}$(p,q)$ against \textsc{arma}$(p+r,q+s)$.

Note that if $s\neq t$, $U_{t}$ and $U_{s}$ are independent (under the null of no serial dependence in the latent process). Also $E(U_{t})=0$ and hence $E(S_{\psi_a}(\beta,\tau,0))=0$. Thus under $H_{0}^{(\psi)}$, $\theta_{0}=(\beta_0,\tau_0,0)$, the $S_{\psi_a}$ are uncorrelated with covariance%
\begin{equation}\label{eqn: Omega aa}
\Omega_{aa} =\mathrm{Var}(S_{\psi_a}(\beta,\tau, 0))=\tau ^{2}\sum_{t=1}^{n}
E\left( U_{t}^{2}\right)  E\left( U_{t-a}^{2}\right).
\end{equation}
The $E\left( U_{t}^{2}\right)$ is calculated using the marginal evaluations
\begin{equation}\label{eq: E Ut2}%
E(U_{t}^2) = \sum_{y_t=0}^{m_t} f(y_{t};\theta_{0}) u_{t}^2,
\end{equation}
where
\begin{equation}\label{eq: YtDensity}
f(y_{t};\beta,\tau) = \int_{\mathbb{R}} f(y_{t}|W_{t}(\beta,\tau)) g(\tilde\alpha_{t})d\tilde\alpha_{t},
\end{equation}
\begin{equation}\label{eq: UtMarginal}%
u_{t}= y_{t}- f^{-1}(y_{t};\beta,\tau)\int_{\mathbb{R}}m_{t} \dot{b}(x_{t}^{\T}\beta+\sqrt{\tau}\tilde\alpha_{t})f(y_{t}|W_{t}(\beta,\tau))
g(\tilde\alpha_{t})d\tilde\alpha_{t}
\end{equation}
and $f(y_{t}|W_{t})$ is given in \eqref{eq: EFDensity}. The integrals are calculated numerically with the \textbf{R}-function ``\textsf{integrate}".

The score statistic is then constructed as
\[
\hat Q_{\psi}(L)=\sum_{a=1}^{L} \frac{n^{-1}S_{\psi_a}^{2}(\hat{\beta}^{(1)},\hat{\tau}^{(1)},0)}{n^{-1}\hat{\Omega}_{aa}}%
\]
where $\hat{\Omega}_{aa}$ is computed with $(\hat\beta^{(1)}, \hat\tau^{(1)}, 0)$ replacing $\theta$ using \eqref{eq: E Ut2}, \eqref{eq: YtDensity} and \eqref{eq: UtMarginal}.

The asymptotic distribution of $\hat Q_{\psi}(L)$ is given in Theorem \ref{Thm: Nulldist existence depend} below. For this it is required that the marginal estimates $(\hat \beta^{(1)}, \hat \tau ^{(1)})$ are consistent and asymptotically normally distributed under the null hypothesis $H_0^{(\psi)}: \psi=0$.
In \cite{DunsHe2016} it is shown that $n^{-1}l_1(\delta)\to Q(\delta)$, where $Q(\delta)$ is defined for both types of regressors. For asymptotic identifiability we require
\begin{condition}\label{Cond: Ident and Const}
$Q(\delta)$ has a unique maximum at the true value $\delta_{0}$.
\end{condition}
This condition holds for most regression sequences described by Condition \ref{Cond: Xt} and sequences of trials satisfying Condition \ref{Cond: mt} -- see \cite{DunsHe2016}.
As shown in \cite{DunsHe2016}, the marginal estimates are consistent and asymptotically normally distributed even when the latent process has serial dependence with the requirement that $\alpha_{t}$ be stationary and strongly mixing with a mixing rate that converges to zero sufficiently rapidly. This is satisfied for the ARMA models in particular. For our present purposes, we only need the asymptotic properties under the hypothesis of independent $\tilde \alpha_t$.

\begin{theorem}\label{Thm: Nulldist existence depend}%
Assume Conditions 1 to 4 and $\tau >0$ and $H_{0}^{(\psi)}: \psi=0$ is true. Then, as $n\to\infty$,
$\hat Q_{\psi}(L)\rightarrow\chi^{2}(L)$ in distribution.
\end{theorem}

Under the null hypothesis of $\psi=0$, $S_{\psi_a}(\beta,\tau)$ and $S_{\psi_b}(\beta,\tau)$ are independent for any $a\neq b$, $a,b=1,\ldots,L$. By establishing the conditions required in Theorem 27.2 in \cite{billingsley1968convergence} it can be shown that $n^{-1/2}S_{\psi_a}(\delta_0)$ is asymptotically normally distributed with mean zero and finite covariance. Next, note that based on the ARMA\ model for $\alpha_t$, it is straightforward to show $E(\tilde\alpha_t|y_t)$ is uniformly bounded, so $\hat\delta^{(1)}\to \delta_{0}$ in probability, and hence using a Taylor expansion, it can be shown that $n^{-1/2}\left(S_{\psi_a}(\hat\delta^{(1)})- S_{\psi_a}(\delta_0)\right)=o_{p}(1)$.

A combined test of $H_{0}:\tau=0,\psi=0$ is not feasible, as the first derivative and the second derivative of the log-likelihood with respect to $\psi$ are both zero vectors when $\tau=0$ so they cannot be used.

\section{Simulation} \label{Sec: Simulations}

In this section we present some simulation results to illustrate the accuracy of asymptotic distributions of score tests for finite samples. Throughout this section the state equation is the linear trend in time with the latent process
\begin{equation}\label{eq: Wt}%
W_{0,t}= 1 + 2(t/n) + \alpha_t
\end{equation}
where $\alpha_t = \phi \alpha_{t-1} + \epsilon_t$, $\epsilon_{t}\sim N(0,\sigma^2_{\epsilon})$, and $\tau=\mathrm{Var}(\alpha_t)$.

For each experiment 10,000 replications were used. For each replication the latent process $\{\alpha_t\}$ is simulated and the observations are generated by
\begin{equation}\label{eq: Yt|alpha}%
Y_t|\alpha_t \sim B(m_t, \pi_{t});\quad \pi_{t}=1/(1+\exp(-W_{0,t})).
\end{equation}
The marginal likelihood estimates of $\hat\beta^{(1)}$ and the square root of $\hat\tau^{(1)}$ were obtained using the \textbf{R} package ``\textsf{lme4}".

The first simulation compares the finite sample distributions of the supremum score statistic $Q_{\tau}(\Psi)$ in \eqref{eq: sup ST} with the theoretical distribution $\mathcal{F}_{\Psi}(u)$ in \eqref{eq: Sup S1 dist}, under the null hypothesis of $\tau=0$. In each replication, the samples of 
$\{Y_t: 1\le t\le n\}$ are generated from model \eqref{eq: Yt|alpha} with $\alpha_t=0$. The supremum score statistics (theoretical or empirical) are obtained over the discrete grid $\Psi=-0.9(0.1)0.9$. Table \ref{tb: quantiles of sup chi2} shows that the empirical quantiles for 10\% and 5\% levels are in good agreement with the theoretical quantiles and at the $2.5\%$ and $1\%$ levels, the theoretical quantiles underestimate the empirical quantiles for $n=200$ sample size.
\begin{table}[ptb]\centering
\caption{Quantiles of theoretical supremum $\chi^2(1)$ distribution and empirical
$\hat Q_{\tau}(\Psi)$ over the scale of $\Psi=[-0.9,0.9]$.}
\begin{tabular}{l*{10}{c}r}
&&\multicolumn{4}{c}{$m_t=1$} & \multicolumn{4}{c}{$m_t=2$}\\
& & 10\% & 5\% & 2.5\% & 1\% & 10\% & 5\% & 2.5\% & 1\%\\
\multirow{2}{*}{$n=200$} & $\mathcal{F}_{\Psi}(u)$& 5.94 & 7.30 & 8.67& 10.48 & 5.04& 6.38 &7.74& 9.53\\
& $\hat Q_{\tau}(\Psi)$ & 5.55 & 7.12 & 9.08 & 12.15 &5.06 & 6.63 & 8.69& 11.86\\
\multirow{2}{*}{$n=10^3$} & $\mathcal{F}_{\Psi}(u)$ & 6.03 &7.39 &8.76 &10.57 &5.43 &6.78 &8.14 & 9.94\\
& $\hat Q_{\tau}(\Psi)$ & 5.56 & 7.07 & 8.48 &10.27&5.26& 6.75& 8.25 & 10.66  \\
\end{tabular}
\label{tb: quantiles of sup chi2}%
\end{table}

\begin{figure}[ptb]\centering
	\subfloat[Binary $(m_{t}=1)$]{\includegraphics[trim=0 0 0 50,clip,width=0.5\textwidth] {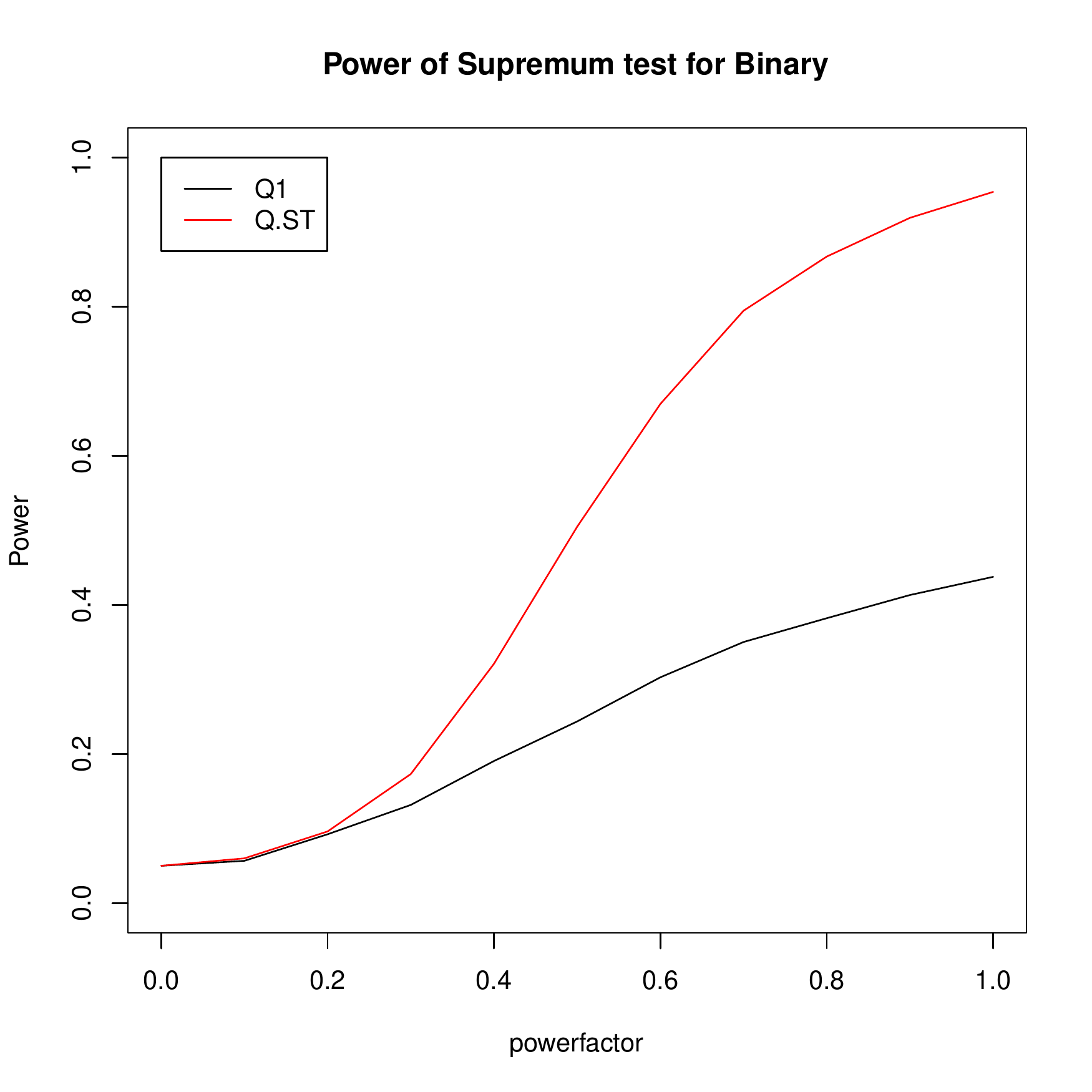}}
	\subfloat[Binomial $(m_{t}=2)$]{\includegraphics[trim=0 0 0 50,clip,width=0.5\textwidth] {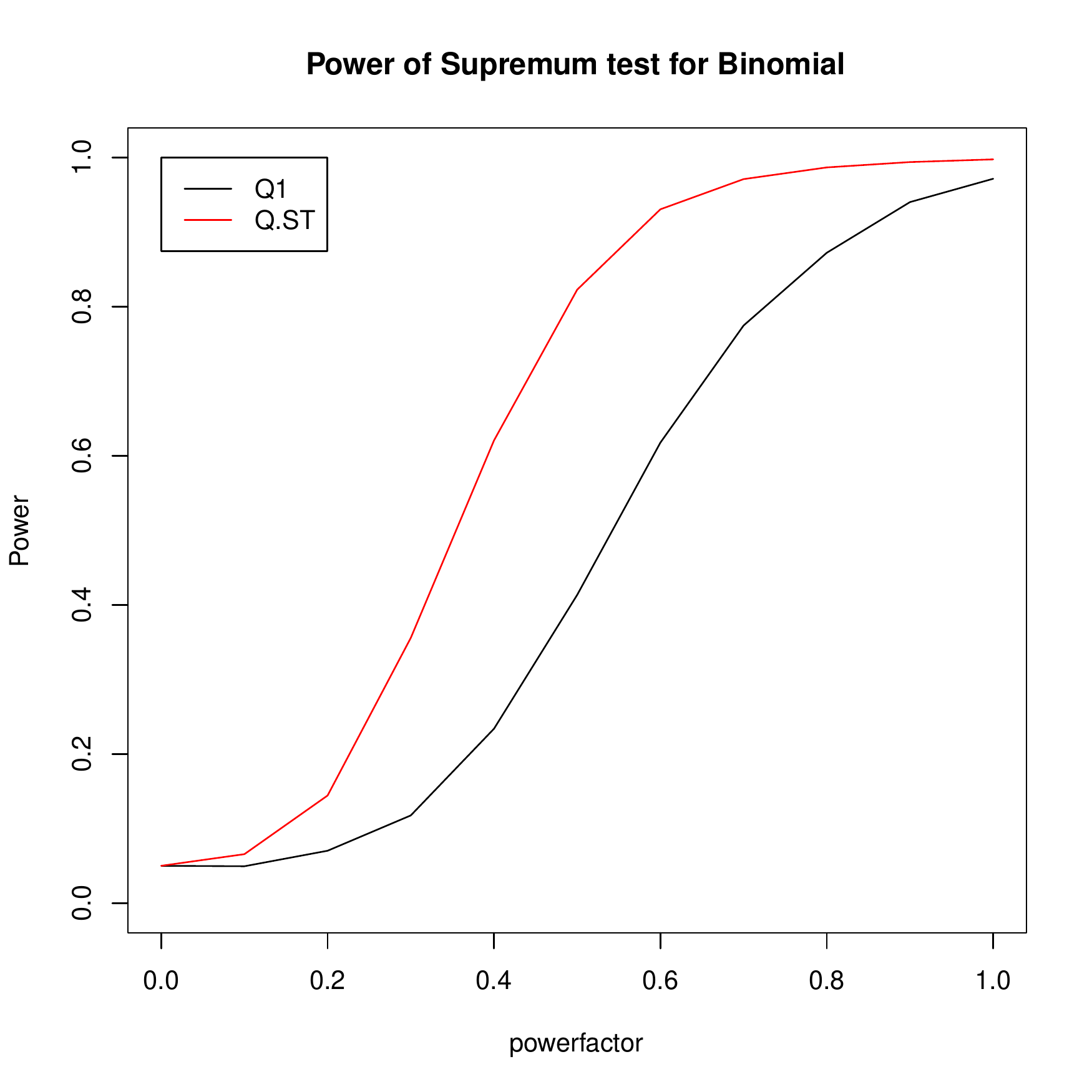}}
	\caption{Power comparison of $\hat Q_{\tau}(\Psi)$ ($Q.ST$ in red) and $\hat Q_{\tau}(0)$ ($Q1$ in black) with Binary(left) and Binomial(right) responses, where  $\Psi=-0.9(0.1)0.9$. The ``powerfactor" refers to value $\sqrt{\tau} = \mathbf{i}\sqrt{\tau_0}$ where $\mathbf{i}=0(0.1)1$. }
\label{Fg: Power PDsup}%
\end{figure}

Next we compare the power of the supremum score test $\hat Q_{\tau}(\Psi)$ with that of the ``standard" score test $\hat Q_{\tau}(0)$ in \eqref{eq: std ST}. Again, $\Psi=-0.9(0.1)0.9$. Each test statistic is simulated with $n=200$, the power is evaluated with the empirical probability that the test statistic exceeds the empirical $95\%$ quantile of the null distribution, under which $\tau=0$. The power for the two tests are evaluated at an increasingly more distant set of alternatives $H_a: \sqrt{\tau} = \mathbf{i}\sqrt{\tau_0}$ where $\tau_0 = 1$ and $\mathbf{i}=0(0.1)1$, in each level of $\mathbf{i}$, $\alpha_t=\tau^{1/2}\tilde\alpha_t$ and $\tilde\alpha_t=0.9\tilde\alpha_{t-1}+\epsilon_t$; $\epsilon_t\sim N(0,1)$. Figure \ref{Fg: Power PDsup} shows that for binary series the power of the ``standard" score test does not increase to 1 as $\tau$ grows. This is because, for binary data, the ``standard" score vector under the alternative $\tau>0$ has an asymptotic normal distribution of mean approximately zero, which is close to its counterpart under the null $\tau=0$, and consequently, the probability to reject the null is small -- see \cite{DunsHe2016}. However, the power improves significantly in the supremum test for both binary and binomial cases. Therefore, the supremum test for serial dependence is recommended particularly for binary series.

The accuracy of the asymptotic distribution of the score test $\hat Q_{\psi}(L)$ for serial dependence will be assessed in one of the examples below. Due to the substantial probability of $\hat\tau^{(1)}=0$ with marginal likelihood estimation, in particular for binary series -- see \cite{DunsHe2016}, the convergence of $\hat Q_{\psi}(L)$ to a chi-squared distribution with binary data is slow.

\section{Applications}\label{Sec: Binary/Binom TS Examples}

In this section we detect, firstly, the existence of a latent process, and if present, the serial dependence, for some real examples, with the methods proposed in this paper. In the test for
a latent process, the alternative $\alpha_{t}$ is assumed to be an AR$(1)$ process: $\alpha_{t}=\phi
\alpha_{t-1} +\varepsilon_{t}$, $\varepsilon_{t}\sim N(0,1)$, where $\psi=\phi$ is the nuisance parameter. Throughout this section we simulate with 10,000 replicates. The supremum score tests
again use $\Psi=-0.9(0.1)0.9$.

\vspace{2mm}

\noindent\centerline{\textbf{\textit{Example 1: Oxford-Cambridge Boat race -- Binary Series}}}

\vspace{1mm}

\cite{klingenberg2008regression} consider the time series of 153 observations over the period 1829 to 2007 of outcomes of the Cambridge-Oxford annual boat race with $y_{t}= 1$ when Cambridge wins and $y_{t} = 0$ otherwise. \cite{klingenberg2008regression} fits a parameter driven regression consisting of an intercept and the single covariate $x_{t}$ being the weight difference between the winning and losing side with an \textsc{ar}$(1)$ latent process. His method allows for time gaps, most of which occur early in the series.  His fitted model implies the presence of substantial serial dependence and so we use this series as a way of illustrating the performance of the statistics defined in this paper. However, since the above tests require equal time spacing, for this application, time is taken to be the sequence number of each race.

We first implement the score test for a latent process under the null $H_{0}:\tau=0$. In each simulation, the binary series is generated with the probability of success: $\pi_{t}=1/(1+\exp(-x_{t}^{\T}\hat\beta^{(0)}))$, where $\hat \beta^{(0)}= (0.194, 0.118)$ are GLM estimates. Table \ref{tb: Nulldist OxCam Corr} summarises the simulated distributions of the ``standard" score test $\hat Q_{\tau}(0)$ in \eqref{eq: std ST}, the supremum score test $\hat Q_{\tau}(\Psi)$ in \eqref{eq: sup ST}, and their reference distributions $\chi^2(1)$ and $\mathcal{F}_{\Psi}(u)$ in
\eqref{eq: Sup S1 dist}. The table shows that for both the ``standard" score test and the supremum score test, there is upward bias for the 1\% quantiles and downward bias for the 20\%, 10\% and 5\% quantiles. The ``standard" score test is insignificant at 5\% level. The observed value of the supremum test statistic is significant at the 1\% level using either the simulated quantile or the theoretical upper bound quantile. The second test statistic for serial dependence requires the marginal fit, $\hat \delta^{(1)}$. For the boat race series, $\hat\tau^{(1)}=0$, hence the test for serial dependence cannot be constructed. As is noted in \cite{DunsHe2016}, the marginal likelihood estimates of binary data can be misleading because the `pile-up' effect happens with approximately a 50\% of chance.

\begin{table}[ptb]\centering
\caption{Null distribution quantiles of score tests from the Cambridge-Oxford boat race series}
\begin{tabular}{l*{7}{c}r}
\text{Test} & & \multicolumn{4}{c}{Distribution} & \multicolumn{1}{c}{Observed}\\
& & 20\% & 10\% & 5\% &1\% & \\
\multirow{2}{*}{Standard} & $\chi^2(1)$  & 1.64& 2.71&  3.84& 6.63  &- \\
& $\hat Q_{\tau}(0)$ & 1.18 & 1.93 & 3.21 & 9.28 & 0.39\\
\multirow{2}{*}{Supremum} & $\mathcal{F}_{\Psi}(u)$ &4.66&6.02& 7.38&10.57 & - \\
& $\hat Q_{\tau}(\Psi)$ & 3.38 & 4.86& 6.70& 11.61& 13.40$^\ast$\\
\end{tabular}
\label{tb: Nulldist OxCam Corr}%
\end{table}

\vspace{2mm}

\noindent \centerline{\textbf{\textit{Example 2: Crime Records -- Binomial Time Series}}}

\vspace{1mm}

\cite{dunsmuir2008assessing} considered the number of convictions, $y_{t}$, obtained from monthly numbers of trials, $m_{t}$, in the higher court in the Australian state of New South Wales, for 6 crime categories: Assault, Sexual Assault, Robbery, Break and Enter, Motor Theft and Other Theft for the
period Jan, 1995 to Jun, 2007. For these series the binomial distribution for the number of charges which led to a successful prosecution is used. For each crime, the regressors $X_{t}=(1, T_{t}, \texttt{DNA}_{t-L},\texttt{SD}_t)$ are defined as: $T_{t}=t/12$ where $t$ is the month since Jan, 1995; $\texttt{DNA}_{t-L}=\max(t-L-73, 0)$ is a linear effect since Jan, 2001 ($t=73$) mirroring the nearly linear growth in the number of individuals in the DNA database, and $L$ is the delay effect of each crime; $\texttt{SD}_t$ represents any seasonal dummy variables.

Table \ref{tb: Nulldist HC convict latent} explores the existence of an latent process for each crime category with the supremum score test $\hat Q_{\tau}(\Psi)$. The simulated null distribution quantiles of the supremum score test are given against which the observed statistics can be compared. The latent process is detected in all crimes except for Motor Theft. Therefore the further test for serial dependence,
$\hat Q_{\psi}(L)$, is justified. $L=2$ was selected because the residuals from the GLM fit suggested at most 2 lags were needed. The independent binomial samples are generated with $Y_{t}|\alpha_{t}\sim B(m_{t}, 1/(1+\exp(-x_{t}^{\T}\hat \beta^{(0)}-\alpha_{t}))$, $\alpha_{t}\sim N(0,1)$, using the GLM fit $\hat\beta^{(0)}$. The testing results show that Break and Enter, Robbery exhibit significant serial dependence.
\begin{table}[ptb]\centering
\caption{Simulated null distribution of the two-step score test from the Higher Court convictions of New South Wales, Australia \\ ( $^\ast$ significance at the $5\%$ level).}
\begin{tabular}{l*{12}{c}r}
& \multicolumn{5}{c}{Test for Latent Process} & \multicolumn{5}{c}{Test for Serial Dependence} \\
& 20\% & 10\%  &  5\%  & 1\%  & \text{observed} & 20\% & 10\%   &   5\%    & 1\%  & \text{observed}\\
\text{Assault} &3.15 & 4.32& 5.68 & 8.79 & 14.81$^\ast$ & 3.33& 4.64 & 6.32 & 9.52 & 0.62 \\
\text{SexAssault} & 3.02 & 4.12& 5.32& 9.33 &91.92$^\ast$ & 3.34 & 4.79 & 6.29& 10.91& 5.05\\
\text{BreakEnter} & 3.40& 4.61& 5.73 & 9.26 &68.99$^\ast$ & 3.19 & 4.69 &6.35 & 9.92&11.56$^\ast$\\
\text{Robbery} & 3.10& 4.29& 5.74 & 9.46& 90.65$^\ast$ & 3.35& 4.79 & 6.17 & 9.52& 8.67$^\ast$\\
\text{MotorTheft}& 3.31 & 4.48 & 5.45 & 9.87& 3.48 & - & - & - & - & -\\
\text{OtherTheft}& 3.09 & 4.44 & 6.11 & 9.68 & 15.33$^\ast$& 3.42 & 5.07& 6.27& 10.10 & 3.32\\
\end{tabular}
\label{tb: Nulldist HC convict latent}%
\end{table}

\section{Conclusions} \label{Sec: Conclusions}

For parameter driven models, we have proposed a pair of score-type tests for, first, the detection of a latent process and, second, serial dependence within it. In the first step the correlation coefficients $\psi$ of latent process are not estimable under the null hypothesis of no latent process and thus are nuisance parameters. Two ways are proposed to deal with this issue: set $\psi=0$ or use a supremum test. The former results in a standard score test evaluated with GLM estimates only, the latter requires the maximum value of score test statistics over the space of nuisance parameters. Simulations show that for binary data, the score test obtained by setting $\psi = 0$ is underpowered for alternative hypothesis. Therefore a supremum score test is necessary when detecting the existence of a latent process in a binary sequence and is capable of providing sufficient power under the alternatives. The simulated quantiles and the theoretical quantiles obtained under the numerical approximation \eqref{eq: Sup S1 dist} for the supremum score test are in good agreement. In practice, the theoretical upper bound quantiles of \eqref{eq: Sup S1 dist} are evaluated with the estimate of $\beta_{0}$.

The score test for serial dependence is constructed against the alternative that the latent process
follows an \textsc{arma}$(p,q)$ process. To establish the asymptotic distribution of the score test for serial dependence, the asymptotic normality of marginal likelihood estimators is required. \cite{DunsHe2016} shows that the marginal likelihood estimators are unbiased and asymptotically normal even if the latent process is correlated, from which it follows that the score statistic $\hat Q_{\psi}(L)$ has an asymptotic chi-squared distribution under the null of no serial dependence.

\section{Appendix}

\centerline{\textit{Outline proof of Theorem \ref{Thm: Nulldist existence latent}}}

Given that $\hat\beta^{(0)}-\beta_{0}=o_{p}(1)$, using a Taylor series expansion it can be shown that
\[
n^{-1/2}\left(S_{\tau,1}(\hat\beta^{(0)})- S_{\tau,1}(\beta_{0})\right) - \sqrt{n}(\hat\beta^{(0)} -\beta_{0})^{T}J_{n}(\beta_0) =o_{p}(1)
\]
in which $n^{-1/2}(\hat\beta^{(0)} - \beta_{0}) - I_{n}^{-1}(\beta_0)U_{n}(\beta_0) =o_{p}(1)$. Note also
\begin{equation*}
n^{-1/2}S_{\tau,2}(\beta,\psi)= n^{-1/2}\sum_{t=2}^{n} e_t(\beta,0) \sum_{h=1}^{t-1}R(h;\psi)e_{t-h}(\beta,0),
\end{equation*}
and $e_t(\beta,0)$ are uniformly bounded for all $\beta$, for any fixed $\psi$,
\begin{equation*}
n^{-1/2}\left( S_{\tau,2}(\hat\beta^{(0)},\psi)- S_{\tau,2}(\beta_0,\psi) \right)=o_{p}(1).
\end{equation*}
Then asymptotically $n^{-1/2}(S_{\tau}(\hat\theta^{(0)}) - S_{\tau}^{\dag}(\theta_0))=o_{p}(1)$.

Observe that $ n^{-1/2}S_{\tau}^{\dag}(\theta_0)$ can be represented as $\sum_{t=1}^n\xi_{nt}$ which is a sum across rows in a  triangular array of martingale differences. Using the central limit theorem for such arrays as in \cite{hall1980martingale} gives
\begin{equation*}
\sum_{t=1}^n \xi_{nt}\overset{d}\to N(0,\underset{n\to\infty}\lim V_{n}(\beta_0,\psi)).
\end{equation*}

\bibliography{TestSerialDepParDrivenBinTSarXiv}
\bibliographystyle{imsart-nameyear}

\end{document}